\documentclass[12pt]{amsart}
\usepackage{amssymb}

\begin{document}


{\theoremstyle{plain}%

\newtheorem{theorem}{Theorem}[section]
  \newtheorem{corollary}[theorem]{Corollary}
  \newtheorem{proposition}[theorem]{Proposition}
  \newtheorem{lemma}[theorem]{Lemma}
  \newtheorem{question}[theorem]{Question}

  \newtheorem{conjecture}[theorem]{Conjecture}
  \newtheorem{claim}[theorem]{Claim}
}
{\theoremstyle{remark}
  \newtheorem{fact}{Fact}
}
{\theoremstyle{definition}
\newtheorem{definition}[theorem]{Definition}
\newtheorem{remark}[theorem]{Remark}
\newtheorem{example}[theorem]{Example}
}

\newcommand{\Tor}{\operatorname{Tor}}
\newcommand{\reg}{\operatorname{res-reg}}
\newcommand{\Reg}{\operatorname{\bf res-reg}}
\newcommand{\nreg}{\operatorname{reg}}
\newcommand{\nReg}{\operatorname{\bf reg}}
\renewcommand{\a}[1]{\operatorname{end}_{#1}}
\newcommand{\bfa}[1]{\frak{a}^{[#1]}}
\newcommand{\brak}[1]{{[#1]}}
\newcommand{\crak}[1]{{\{#1\}}}
\newcommand{\beg}{\operatorname{beg}}
\newcommand{\ann}{\operatorname{ann}}
\newcommand{\ass}{\operatorname{Ass}}

\def\N{\mathbb{N}}
\def\Z{\mathbb{Z}}
\def\F{\mathbb{F}}
\def\GG{\mathbb{G}}
\def\M{\mathcal{M}}
\def\K{\mathbb{K}}
\def\R{\mathcal{R}}
\def\X{\mathbf{X}}
\def\Y{\mathbf{Y}}
\def\V{\mathbb{V}}
\def\kk{\mathbf{k}}
\def\b{\mathbf{b}}
\def\x{\mathbf{x}}
\def\y{\mathbf{y}}
\def\z{\mathbf{z}}
\def\e{\mathbf{e}}
\def\f{\mathbf{f}}
\def\0{\mathbf{0}}
\def\1{\mathbf{1}}
\def\n{\mathbf{n}}
\def\m{\mathbf{m}}
\def\k{\mathbf{k}}
\def\d{\mathbf{d}}
\def\r{\mathbf{r}}
\def\u{\mathbf{u}}

\newcommand{\msg}[1]{\marginpar{\addtolength{\baselineskip}{-3pt}{\footnote size \it #1}}}
\numberwithin{equation}{section}


\title[Asymptotic behavior of multigraded regularity]{Minimal free resolutions and asymptotic behavior of multigraded regularity}
\author{Huy T\`ai H\`a}
\address{Tulane University \\
Department of Mathematics \\
6823 St. Charles Avenue \\
New Orleans, LA 70118, USA} \email{tha@tulane.edu}
\urladdr{http://www.math.tulane.edu/~tai}
\thanks{The first author is partially supported by Louisiana Board of Regents Enhancement Grant.}

\author{Brent Strunk}
\address{University Of Louisiana at Monroe \\
Department of Mathematics and Physics \\
700 University Avenue \\
Monroe, LA 71209, USA} \email{bstrunk@ulm.edu}
\urladdr{http://www.math.purdue.edu/~bstrunk}


\begin{abstract}
Let $S$ be a standard $\N^k$-graded polynomial ring over a field $\kk$, let
$I$ be a multigraded homogeneous ideal of $S$, and let $M$ be a
finitely generated $\Z^k$-graded $S$-module. We prove that the
resolution regularity, a multigraded variant of Castelnuovo-Mumford
regularity, of $I^nM$ is asymptotically a linear function. This
shows that the well known $\Z$-graded phenomenon carries to the
multigraded situation.
\end{abstract}

\maketitle


\section{Introduction}

In this paper we investigate the asymptotic behavior of a
multigraded variant of Castelnuovo-Mumford regularity, the
resolution regularity.

Throughout the paper, $\N$ will denote the set of non-negative integers.
Let $\e_1, \dots, \e_k$ be the standard unit vectors of $\Z^k$ and
let $\0 = (0,\dots,0)$. Let $\kk$ be a field and let $S$ be a
standard $\N^k$-graded polynomial ring over $\kk$. That is, $S_\0 = \kk$ and
$S$ is generated over $ S_\0$ by elements of $\bigoplus_{l=1}^k
S_{\e_l}$. For $l=1, \dots, k$, suppose $S_{\e_l}$ is generated as a
vector space over $S_\0$ by $\{x_{l,1 }, \dots, x_{l,N_l}\}$. For
simplicity, we shall use $\x_l$ to represent the elements
$\{x_{l,1}, \dots, x_{l,N_l}\}$ and $(\x_l)$ to denote the $S$-ideal
generated by $\{x_{l,1}, \dots, x_{l,N_l}\}$.

\begin{definition} \label{def.reg}
Let $M = \bigoplus_{\n \in \Z^k}M_\n$ be a finitely generated
$\Z^k$-graded $S$-module. Let $$\F: 0 \rightarrow F_p \rightarrow
\dots \rightarrow F_1 \rightarrow F_0 \rightarrow M \rightarrow 0$$
be a minimal $\Z^k$-graded free resolution of $M$ over $S$, where
$F_i = \bigoplus_j S(-c^\brak{1}_{ij}, \dots, -c^\brak{k}_{ij})$ for
$i=1, \dots, p$. For each $1 \le l \le k$, set $c^\brak{l}_i =
\max_j \{c^\brak{l}_{ij}\}$ for $i=0, \dots, p$ and let
$$\reg_l(M) = \max_i \{c^\brak{l}_i - i\}.$$
The multigraded {\it resolution regularity} of $M$ is defined to be
the vector
$$\Reg(M) = (\reg_1(M), \dots, \reg_k(M)) \in \Z^k.$$
\end{definition}

\begin{example} Let $S = \mathbb{Q}[x_0,x_1,y_0,y_1]$, where
$\deg x_i = \e_1 \in \Z^2$ and $\deg y_j = \e_2 \in \Z^2$ for all $0
\le i,j \le 1$. Then $S$ is a standard $\N^2$-graded polynomial ring over
$\mathbb{Q}$. Consider the ideal $I=(x_0^2, x_0y_1, x_1y_0, y_0^2)$
as a $\Z^2$-graded $S$-module. Using any computational algebra
package (such as, CoCoA \cite{cocoa} or Macaulay2 \cite{Mac2}), we
obtain the following minimal $\mathbb{Z}^2$-graded free resolution
of $I$ over $S$:
$$0\rightarrow S(-3,-3)\rightarrow \begin{array}{c}
S(-3,-2)^2\\
\bigoplus\\
S(-2,-3)^2\\
\end{array} \rightarrow \begin{array}{c}
S(-1,-2)\\
\bigoplus\\
S(-2,-1)\\
\bigoplus\\
S(-3,-1)\\
\bigoplus\\
S(-2,-2)^2\\
\bigoplus\\
S(-1,-3)\\
\end{array} \rightarrow \begin{array}{c}
S(-2,0)\\
\bigoplus\\
S(-1,-1)^2\\
\bigoplus\\
S(0,-2)\\
\end{array}\rightarrow I\rightarrow 0.$$
By definition, $\reg_1(I) = 2$ and $\reg_2(I) = 2$. Thus,
$\Reg(I)=(2,2) \in \Z^2$.
\end{example}

\begin{remark} \label{tor}
For a vector $\n \in \Z^k$ we shall use $n_l$ to denote its $l$-th coordinate. It can be seen that $\Tor^S_i(M,\kk) = H_i(\F \otimes \kk)$. Thus,
the resolution regularity of $M$ can also be calculated as follows.
For each $1 \le l \le k$,
$$\reg_l(M) = \max \{ n_l ~|~ \exists \n \in \Z^k, i \ge 0 \text{ so that }
\Tor^S_i(M,\kk)_{\n+i\e_l} \not= 0\}.$$
We shall often make use of this observation.
\end{remark}

The resolution regularity is a multigraded variant of the well known
Castelnuovo-Mumford regularity, developed via the theory of Hilbert
functions and minimal free resolutions by Aramova, Crona and DeNegre
\cite{acn} (for $\Z^2$-gradings), and Sidman and Van Tuyl \cite{sv}
(for $\Z^k$-gradings in general). This complements the notion of {\it
multigraded regularity}, another variant of Castelnuovo-Mumford
regularity, studied by Hoffman and Wang \cite{hh} (for
$\Z^2$-graded), and Maclagan and Smith \cite{ms1} (for $G$-graded,
where $G$ is an Abelian group).
Roughly speaking, the resolution regularity of $M$ captures the
maximal coordinates of minimal generating multidegrees of syzygy
modules of $M$. In particular, it provides a crude bound for the
generating degrees of $M$. Thus, resolution regularity can be viewed
as a refinement of Castelnuovo-Mumford regularity. Resolution
regularity, furthermore, shares important similarities with the
original definition of Castelnuovo-Mumford regularity (cf.
\cite{ha3}).

Many authors \cite{acn,hv,ha3,hh,ms1,ms2,romer,sv,SVW} in recent
years have tried to extend our knowledge of Castelnuovo-Mumford
regularity to the multigraded situation. On the other hand, there has
been a surge of interest in the asymptotic behavior of
Castelnuovo-Mumford regularity of powers of an ideal in an
$\N$-graded algebra (cf. \cite{BEL, Chandler, cht, Cut, CuEL, ggp,
K, Swanson, TW}). It is known that if $S$ is a standard $\N$-graded
algebra (over a Noetherian ring $A$), $I$ is a homogeneous ideal in
$S$ and $M$ is a finitely generated $\Z$-graded $S$-module, then
$\operatorname{reg} (I^nM)$ is asymptotically a linear function in
$n$ with slope $\le d(I)$, where $d(I)$ is the maximal generating
degree of $I$. The aim of this paper is to show that the $\Z$-graded
phenomenon carries to the multigraded setting.

Suppose $M$ is a $\Z^k$-graded $S$-module, minimally generated in degrees $\d_1(M), \dots, \d_v(M)$, where $\d_i(M) = (d_{i,1}, \dots, d_{i,k})$. Then, for each $1 \le l \le k$, we define
$$d^\brak{l}(M) = \max \{d_{i,l} ~|~ i = 1, \dots, v\}$$
to be the maximal $l$-th coordinate among the minimal generating degrees of $M$. Our first main result is stated as follows.

\begin{theorem}[Theorem~\ref{asymptotic linearity}] \label{intro.thm1}
Let $S$ be a standard $\N^k$-graded
polynomial ring over a field $\kk$ and let $M$ be a
finitely generated $\Z^k$-graded $S$-module. Let $I$ be a
multigraded homogeneous ideal in $S$ minimally generated in degrees $\d_1(I),
\dots, \d_v(I)$. Then, $\Reg(I^nM)$ is asymptotically a linear
function with slope vector at most
$\big( d^\brak{1}(I), \dots, d^\brak{k}(I) \big)$ componentwise.
\end{theorem}

The slope vector of $\Reg(I^nM)$, in fact, can be described explicitly via the theory of {\it reductions}.

\begin{definition} \label{def.reduction}
We say that a multigraded homogeneous ideal $J \subset I$ is an {\it $M$-reduction} of $I$ if $I^nM = JI^{n-1}M$ for all $n \gg 0$. For each $l = 1, \dots, k$, we define
$$\rho^\brak{l}_M(I) = \min \{ d^\brak{l}(J) ~\big|~ J \mbox{ is an } M-\mbox{reduction of } I\}$$
(here, $d^\brak{l}(J)$ is defined similarly to $d^\brak{l}(I)$), and
$$\beg^\brak{l}(M) = \min \{ n_l ~\big|~ \exists \n \in \Z^k \mbox{ such that } M_\n \not= 0\}.$$
\end{definition}

\begin{theorem}[Theorem~\ref{coefficient}] \label{intro.thm2}
Let $S$ be a standard $\N^k$-graded polynomial ring over a field $\kk$ and let $M$ be a finitely generated $\Z^k$-graded $S$-module. Let $I$ be a multigraded homogeneous ideal in $S$. Then, there exists a vector $\b \ge (\beg^\brak{1}(M), \dots, \beg^\brak{k}(M))$ such that for $n \gg 0$,
$$\Reg(I^nM) = n \big(\rho^\brak{1}_M(I), \dots, \rho^\brak{k}_M(I)\big) + \b.$$
\end{theorem}

It is well known that in the $\Z$-graded case, Castelnuovo-Mumford
regularity can be defined either via the shifts of a minimal free
resolution or by local cohomology modules with respect to the
homogeneous irrelevant ideal (cf. \cite{eg}). In the multigraded
situation (i.e. when $k \ge 2$), the two approaches, using minimal
free resolutions and local cohomology with respect to the irrelevant
ideal, to define variants of regularity do not agree and thus give
different invariants (see \cite{ha3} for a discussion on the
relationship between these two multigraded variants of regularity).
The incomparability between these approaches makes it difficult to
generalize $\Z$-graded results to the multigraded setting. Another
conceptual difficulty is the simple fact that $\Z^k$ is not a
totally ordered set. Thus, it is not easy to capture maximal
coordinates of shifts in the minimal free resolution by a single
vector.

Our method in this paper is inspired by a recent work of Trung and
Wang \cite{TW}, where the authors study Castelnuovo-Mumford
regularity via {\it filter-regular sequences} (which originated from
\cite{STC}). To start, we first define a multigraded variant of the
notion of filter-regular sequences.
We then need to establish the
correspondence between resolution regularity and filter-regular
sequences, which is no longer apparent after passing to the multigraded
situation. To achieve this, we develop previous work of R\"omer
\cite{romer} and Trung \cite{RedExp} to fuller generality to link
resolution regularity, filter-regular sequences and local cohomology
with respect to different graded ideals, via the theory of
Koszul homology.  Our main theorems are obtained by
generalizing to the multigraded setting techniques of \cite{TW} in investigating filter-regular
sequences. Furthermore, as a byproduct of
our work, it seems possible to define multigraded regularity for
$\Z^k$-graded modules over a standard $\N^k$-graded algebra over an
arbitrary Noetherian ring $A$ (not necessarily a field), even though,
in this case, a finite minimal free resolution as in
Definition~\ref{def.reg} may not exist. This addresses a question
Bernd Ulrich has asked us.

Our paper is structured as follows. In Section~\ref{s.frs}, we
define a multigraded variant of the notion of filter-regular
sequences (Definition~\ref{defn: filter regular}), and generalize properties of filter-regular sequences from
the $\Z$-graded case to the multigraded situation (Lemmas \ref{primes avoidance} and \ref{flat extension}). In
Section~\ref{s.Koszul}, we establish the correspondence between
resolution regularity, filter-regular sequences and local cohomology
with respect to different graded ideals (Theorem~\ref{reg: cohomology and end}).
Theorem~\ref{reg: cohomology and end} lays the groundwork for the rest of the paper.
In Section~\ref{s.asympreg}, we prove our main theorems (Theorems \ref{asymptotic linearity} and \ref{coefficient}).
In the Appendix, making use of our work in Section~\ref{s.Koszul},
we propose an alternative definition for a multigraded variant of Castelnuovo-Mumford
regularity when there
might not exist finite minimal free resolutions (Definition~\ref{def.alt}).

Now, we shall fix some notations and terminology. If $f \in M$ is a homogeneous element of
multidegree $\d = (d_1, \dots, d_k)$, then we define
$$\deg_l(f) = d_l$$
to be the $l$-th coordinate of the multidegree of $f$.
Since coordinates of multidegrees will be discussed frequently in
the paper, we point out that the operator $\bullet^\brak{l}$
indicates that the maximal $l$-th coordinate of involved multidegrees is
being considered. For $\m = (m_1, \dots, m_k) \in \Z^k$ and $n \in
Z$, by writing $(\m,n)$ we refer to $(m_1, \dots, m_k,n) \in
\Z^{k+1}$.

\vspace{1.5ex}
\noindent{\bf Acknowledgment:} This project started when hurricane Katrina was about to hit New Orleans, where the two authors were at the time. The collaboration was then made possible with the help of many people from various places. We would like to thank the University of Missouri-Columbia and the University of Arkansas for their hospitalities. Especially, we would also like to express our thanks to Ian Aberbach, Jean Chan, Dale Cutkosky, Mark Johnson, Hema Srinivasan and Bernd Ulrich. We thank an anonymous referee for a careful reading of the paper.


\section{Filter-regular sequences} \label{s.frs}

In this section, we extend the notion of filter-regular sequences to the
multigraded situation. More precisely, we define $M$-filter-regular
sequences with respect to a given coordinate.

The first step in our proof of Theorem \ref{intro.thm1} is to express $\reg_l(I^nM)$ in terms of invariants associated to $M$-filter-regular sequences with respect to the $l$-th coordinate. To achieve this we show that, by a flat base extension, $\x_l$ can be taken to be an $M$-filter-regular sequence with respect to the $l$-th
coordinate for all $1 \le l \le k$.

\begin{definition} \label{defn: a-invariant}
Let $M = \bigoplus_{\n \in \Z^k} M_\n$ be a finitely generated
$\Z^k$-graded $S$-module. For each $1 \le l \le k$, define
$$\a{l}(M) = \max \{ n_l ~|~ \exists \n \in \Z^k : M_\n \not= 0\}.$$
\end{definition}

\begin{remark} In general, $\a{l}(M)$ may be infinite, and $\a{l}(M) = - \infty$ if and only if $M = 0$.
\end{remark}

\begin{example} \label{RunExample1}
Consider $S = \kk[x,y,z]$ as a standard $\N^3$-graded polynomial
ring, where $\deg(x) = (1,0,0)$, $\deg(y) = (0,1,0)$, and
$\deg(z) = (0,0,1)$. Let $M = S/I$ where $I=(x^2y,xy^2,xyz, y^3,
y^2z,yz^2)$. Then $\a{1}(M) = \infty$, $\a{2}(M) = 2$ and $\a{3}(M)
= \infty$.
\end{example}

\begin{definition} \label{defn: filter regular}
Fix an integer $1 \le l \le k$. A sequence $\f = \langle f_1, \dots,
f_s \rangle$ of multigraded homogeneous elements in $S$ is said to
be an {\it $M$-filter-regular sequence with respect to the $l$-th coordinate}
if
\begin{enumerate}
\item[$(i)$] $(f_1, \dots, f_s)M \not= M$.
\item[$(ii)$] For every $1 \le i \le s$, we have
$$\a{l}\Big( \dfrac{(f_1, \dots, f_{i-1})M :_M f_i}{(f_1, \dots, f_{i-1})M}
 \Big) < \infty.$$
\end{enumerate}

When $\f = \langle f_1, \dots, f_s \rangle$ is an $M$-filter-regular
sequence with respect to the $l$-th coordinate, we define
$$\bfa{l}_M(\f) = \max_{1 \le i \le s} \Big\{ \a{l}\Big( \dfrac{(f_1, \dots, f_{i-1})M :_M f_i}{(f_1, \dots, f_{i-1})M} \Big) \Big\}.$$
\end{definition}

\begin{example} \label{RunExample2}
Consider $S = \kk[x,y,z]$ and $M = S/(x^2y,xy^2,xyz, y^3,
y^2z,yz^2)$ as in Example \ref{RunExample1}. Then $yM \not= M$ and
$\a{2}(0 :_M y) = 2$. Thus, $y$ is $M$-filter-regular with respect
to the second coordinate. Notice also that $y$ is not
$M$-filter-regular with respect to the first or third coordinate,
while $\langle x,z \rangle$ is an $M$-filter-regular sequence with respect to the
first, second, and third coordinates.
\end{example}

The following lemma gives a characterization for $M$-filter-regular sequences with respect to the $l$-th coordinate via primes avoidance (see also \cite[Lemma 2.1]{RedExp}).

\begin{lemma} \label{primes avoidance}
A sequence $\f = \langle f_1, \dots, f_s \rangle$ of elements in $S$
is an $M$-filter-regular sequence with respect to the $l$-th coordinate if and
only if for each $i = 1, \dots, s$, $f_i \not\in P$ for any
associated prime $P \not\supseteq (\x_l)$ of $M/(f_1, \dots,
f_{i-1})M$.
\end{lemma}

\begin{proof} By replacing $M$ with $M/(f_1, \dots, f_{i-1})M$, it suffices to prove the statement for $i = 1$. For simplicity, we shall write $f$ for $f_1$.

Observe that $\a{l}(0 :_M f) < \infty$ if and only if there exists an integer $m$ such that $(\x_l)^m(0 :_M f) = 0$. In particular, this means that
$$0 :_M f \subseteq \bigcup_{j = 1}^\infty [0 :_M (\x_l)^j] = H^0_{(\x_l)}(M).$$
Let $N = M/H^0_{(\x_l)}(M)$. Observe further that $\ass(M)$ is the disjoint union of $\ass(H^0_{(\x_l)}(M))$ and $\ass(N)$ (cf. \cite[Exercise 2.1.12]{bs}).

Consider an arbitrary $Q \in \ass(M)$ and assume that $Q = \ann_S(h)$ for some $h \in M$. Then, $Q \supseteq (\x_l)$ if and only if $h(\x_l)^j = 0$ for some $j \ge 1$ (since $Q$ is prime), that is, $h \in H^0_{(\x_l)}(M)$. This implies that $Q \in \ass(H^0_{(\x_l)}(M))$. Thus, the set of associated primes of $M$ which do not contain $(\x_l)$ is $\ass(N)$. Hence, it remains to show that $0 :_M f \subseteq H^0_{(\x_l)}(M)$ if and only if $f \not\in P$ for any $P \in \ass(N)$. Indeed, $0 :_M f \subseteq H^0_{(\x_l)}(M)$ if and only if $0 :_N f = 0$ in $N$, that is, $f$ is a non-zerodivisor on $N$, and this is the case if and only if $f$ does not belong to any associated prime of $N$. The lemma is proved.
\end{proof}

The next result allows us to assume that $\x_l$ is $M$-filter-regular with respect to the $l$-th coordinate after a flat extension (see also \cite[Lemma 1.2]{TW}).

\begin{lemma} \label{flat extension}
For $i = 1, \dots, N_l$, let $z_i = \sum_{j=1}^{N_l} u_{ij}x_{l,j}$,
where $U = (u_{ij})_{1 \le i,j \le N_l}$ is a square matrix of
indeterminates. Let $\kk' = \kk(U)$, $S' = S \otimes_\kk \kk'$ and
$M' = M \otimes_\kk \kk'$. Then,
\begin{enumerate}
\item The sequence $\z = \langle z_1, \dots, z_{N_l} \rangle$ is an $M'$-
filter-regular sequence with respect to the $l$-th coordinate. Notice that
$(z_1, \dots, z_{N_l}) = (\x_l)S'$.
\item $\Reg(M) = \Reg(M')$ where $\Reg(M')$ is calculated over $S'$.
\end{enumerate}
\end{lemma}

\begin{proof} To prove (1) we first observe that since the elements $z_1, \dots, z_{N_l}$ are given by independent sets of indeterminates, by
induction it suffices to show that $z_1$ is $M'$-filter-regular
with respect to the $l$-th coordinate. By Lemma \ref{primes avoidance}, we
need to show that $z_1 \not\in P$ for any associated prime $P
\not\supseteq (\x_l)S'$ of $M'$. Indeed, if $P \not\supseteq
(\x_l)S'$ is an associated prime of $M'$, then since $S \rightarrow
S'$ is a flat extension, we must have $P = \wp S'$ for some
associated prime $\wp \not\supseteq (\x_l)S$ of $M$. However, since
$(\x_l)S = (x_{l,1}, \dots, x_{l,N_l})S \not\subseteq \wp$, we must
have $z_1 = u_{1,1}x_{l,1} + \dots + u_{1,N_l} x_{l,N_l} \not\in \wp
S' = P$.

Now, let $\F$ be a minimal $\Z^k$-graded free resolution of $M$ over
$S$. To prove (2), we observe that since $S \rightarrow S'$ is a
flat extension, $\F\otimes_S S'$ is a minimal $\Z^k$-graded free
resolution of $M'$ over $S'$. Thus, by definition, $\Reg(M)
=\Reg(M')$.
\end{proof}


\section{Koszul homology and $a$-invariant} \label{s.Koszul}

The goal of this section is to relate the resolution regularity to
invariants associated to filter-regular sequences. Our techniques
are based upon investigating Koszul complexes and local cohomology
modules with respect to different irrelevant ideals (given by
different sets of variables).

For $\u = (u_1, \dots, u_k) \in \Z^k_{\ge 0}$ such that $u_l \le
N_l$ for all $l = 1, \dots, k$, we shall denote by $\K(\u)$ the
Koszul complex with respect to $(x_{1,1}, \dots, x_{1,u_1}, \dots,
x_{k,1}, \dots, x_{k,u_k})$. Let $H_i(\u;M)$ denote the homology
group $H_i(\K(\u) \otimes_S M)$. When $\u = (N_1, \dots, N_{s-1}, r, 0, \dots, 0)$, for simplicity, we shall write $\K^\crak{s}(r)$ and $H^\crak{s}_i(r;M)$ to denote $\K(\u)$ and
$H_i(\u;M)$, respectively.

\begin{lemma} \label{homology exact sequence}
For each $s$ and $0 \le r \le N_s-1$, let
$$\tilde{H}^\crak{s}_0(r;M) = [0 :_{M/(\x_1, \dots, \x_{s-1}, x_{s,1}, x_
{s,2}, \dots, x_{s,r})M} x_{s,r+1}].$$
Then we have the following exact sequence
\begin{align*}
& \dots \rightarrow H^\crak{s}_i(r;M)(-\e_s)
\stackrel{x_{s,r+1}}{\rightarrow} H^\crak{s}_i(r;M) \rightarrow
H^\crak{s}_i(r+1;M) \rightarrow H^\crak{s}_{i-1}(r;M)(-\e_s) \rightarrow \\
& \dots \rightarrow
H^\crak{s}_1(r;M)(-\e_s)\stackrel{x_{s,r+1}}{\rightarrow}
H^\crak{s}_1(r;M) \rightarrow H^\crak{s}_1(r+1;M) \rightarrow
\tilde{H}^\crak{s}_0(r;M)(-\e_s) \rightarrow 0.
\end{align*}
\end{lemma}

\begin{proof} Consider the following exact sequence of Koszul complexes
\begin{align}
0 \rightarrow \K^\crak{s}(r)(-\e_s)\stackrel{x_{s,r+1}}{\rightarrow}
\K^\crak{s}(r) \rightarrow \K^\crak{s}(r+1) \rightarrow 0.
\label{Koszul sequence}
\end{align}
The conclusion is obtained by taking the long exact sequence of
homology groups associated to (\ref{Koszul sequence}); notice that the last term
$\tilde{H}^\crak{s}_0(r;M)(-\e_s)$ is given by the kernel of the map $H^\crak{s}_0(r;M)(-\e_s) \stackrel{x_{s,r+1}}{\rightarrow} H^\crak{s}_0(r;M)$.
\end{proof}

The following theorem allows us to relate resolution regularity to
invariants associated to filter-regular sequences.

\begin{theorem} \label{reg: filter-regular}
Let $S$ be a standard $\N^k$-graded polynomial ring over a field $\kk$ and
let $M$ be a finitely generated $\Z^k$-graded $S$-module. Suppose
that $\x_l$ is an $M$-filter-regular sequence with respect to the $l$-th
coordinate. Then,
$$\reg_l(M) = \max \{\bfa{l}_M(\x_l), d^\brak{l}(M)\}.$$
\end{theorem}

\begin{proof} Without loss of generality, we shall prove the statement for $l = 1$.

For each
$1\le s \le k$ and $0 \le r \le N_s$, let $Q_{s,r} = N_1 + \dots +
N_{s-1}+r$ and define
$$v^\brak{s}(r) = \max \{ n_1 ~|~ \exists \n \in \Z^k : H^\crak{s}_i(r;M)_{\n+i\e_1} \not= 0 \text{ for some } 1 \le i \le Q_{s,r}\}$$ (we make the convention that $v^\brak{1}(0) = -\infty$).
Observe that $\Tor^S_0(M,\kk)_\d \not= 0$ if and only if $\d$ is a
generating degree of $M$. Thus,
$$d^\brak{1}(M) = \max \{n_1 ~|~ \exists \n \in \Z^k : \Tor^S
_0(M,\kk)_\n \not= 0\}.$$
Now, since $\K^\crak{k}(N_k)$ gives a minimal $\Z^k$-graded free
resolution of $\kk$ over $S$, we have $H^\crak{k}_i(N_k;M) =
\Tor^S_i(M,\kk)$ for all $i \ge 1$. Hence,
$$\reg_1(M) = \max \{v^\brak{k}(N_k), d^\brak{1}(M)\}.$$
It remains to show that
$$\max \{v^\brak{k}(N_k), d^\brak{1}(M)\} =
\max \{ \bfa{1}_M(\x_1), d^\brak{1}(M) \}.$$
To this end, we proceed
in the following steps.
\begin{enumerate}
\item[$(i)$] For any $1 \le r \le N_1$, $$v^\brak{1}(r) = \max_{1 \le i \le r} \Big\{\a{1}\Big( \dfrac{(x_{1,1}, \dots, x_{1,i-1})M :_M x_{1,i}}{(x_{1,1}, \dots, x_{1,i-1})M} \Big) \Big\}.
$$
\item[$(ii)$] For any $s \ge 2$ and any $0 \le r \le N_s$, $$\max\{v^\brak{s}(r), d^\brak{1}(M)\} = \max \{\bfa{1}_M(\x_1), d^\brak{1}(M) \}.$$
\end{enumerate}

For simplicity of notation, we shall denote $\a{1}\Big(
\dfrac{(x_{1,1}, \dots, x_{1,i-1})M :_M x_{1,i}}{(x_{1,1}, \dots,
x_{1,i-1})M} \Big)$ by $s_i$ for $i=1, \dots, N_1$. Then
$\bfa{1}_M(\x_1) = \max \{ s_1, \dots, s_{N_1} \}$.

We shall prove (i) by induction on $r$. By Lemma \ref{homology exact
sequence} we have
$$0 \rightarrow H^\crak{1}_1(1;M) \rightarrow \tilde{H}^\crak{1}_0(0;M)(-\e_1) \rightarrow 0.$$
Thus,
\begin{align*}
v^\brak{1}(1) & = \max \{ n_1 ~|~ \exists \n \in \Z^k:
H^\crak{1}_1(1;M)_{\n+\e_1} \not= 0 \} \\
& = \max \{ n_1 ~|~ \exists \n \in \Z^k:\tilde{H}^\crak{1}_0(0;M)_\n \not= 0\} \\
& = \a{1}(0 :_M x_{1,1}) = s_1.
\end{align*}
This proves the statement for $r=1$.

Suppose that $r > 1$. By Lemma \ref{homology exact sequence}, we
have
$$\dots \rightarrow H^\crak{1}_1(r;M) \rightarrow
\tilde{H}^\crak{1}_0(r-1; M)(-\e_1) \rightarrow 0.$$
This implies that $v^\brak{1}(r) \ge s_r.$ Observe that if
$v^\brak{1}(r-1) = - \infty$ then $v^\brak{1}(r) \ge
v^\brak{1}(r-1)$, which implies that $v^\brak{1}(r) \ge \max
\{v^\brak{1}(r-1), s_r\}.$ Assume that $v^\brak{1}(r-1) > - \infty$.
That is, there exist an integer $i$ and $\n \in
\Z^k$ with $n_1 = v^\brak{1}(r-1)$ such that
$H^\crak{1}(r-1;M)_{\n+i\e_1} \not= 0$ and
$H^\crak{1}(r-1;M)_{\n+(i+1)\e_1} = 0$. By Lemma \ref{homology exact
sequence}, we have
$$\dots \rightarrow H^\crak{1}_{i+1}(r;M)_{\n+(i+1)\e_1} \rightarrow H^\crak{1}_i(r-1;M)_{\n+i\e_1} \rightarrow 0.$$
This implies that
$H^\crak{1}_{i+1}(r;M)_{\n+(i+1)\e_1} \not= 0$. Thus,
$v^\brak{1}(r)\ge v^\brak{1}(r-1)$, and therefore, we also have
$v^\brak{1}(r) \ge \max \{v^\brak{1}(r-1), s_r\}$.

We will prove the other direction, that is, $v^\brak{1}(r) \le \max
\{ v^\brak{1}(r-1), s_r\}$. Consider an arbitrary $\n = (n_1, \dots,
n_k) \in \Z^k$ with $n_1 > \max \{ v^\brak{1}(r-1), s_r \}$. For $i
\ge 2$, by Lemma \ref{homology exact sequence}, we have
$$\dots \rightarrow H^\crak{1}_i(r-1;M)_{\n+i\e_1} \rightarrow H^\crak{1}_i(r;M)_{\n+i\e_1} \rightarrow H^\crak{1}_{i-1}(r-1;M)_{\n+(i-1)\e_1}
\rightarrow \dots$$
Since $n_1 > v^\brak{1}(r-1)$, we have
$H^\crak{1}_{i-1}(r-1;M)_{\n+(i-1)\e_1} =
H^\crak{1}_i(r-1;M)_{\n+i\e_1} = 0$. This implies that
$H^\crak{1}_i(r;M)_{\n+i\e_1} = 0$. For $i=1$, by Lemma
\ref{homology exact sequence}, we have
\begin{align}
\dots \rightarrow H^\crak{1}_1(r-1;M)_{\n+\e_1} \rightarrow
H^\crak{1}_1(r;M)_{\n+\e_1} \rightarrow
\tilde{H}^\crak{1}_0(r-1;M)_\n \rightarrow 0. \label{seqr1}
\end{align}
Since $n_1 > v^\brak{1}(r-1)$, we have
$H^\crak{1}_1(r-1;M)_{\n+\e_1} = 0$. Since $n_1 > s_r$, we have
$\tilde{H}^\crak{1}_0(r-1;M)_\n = 0$. Thus, (\ref{seqr1}) implies
that $H^\crak{1}_1(r;M)_{\n+\e_1} = 0$. Hence, $n_1 >
v^\brak{1}(r)$. This is true for any $n_1 > \max \{v^\brak{1}(r-1),
s_r\}$, so we must have $v^\brak{1}(r) \le \max \{ v^\brak{1}(r-1),
s_r\}.$

We have shown that $v^\brak{1}(r) = \max \{v^\brak{1}(r-1), s_r\}$.
By induction, $$v^\brak{1}(r-1) = \max \{s_1, \dots, s_{r-1}\}.$$
Thus, $v^\brak{1}(r) = \max \{s_1, \dots, s_r\}$ and (i) is proved.

We shall prove (ii) by using double induction on $s$ and on $r$. By
part (i), we have $v^\brak{2}(0) = v^\brak{1}(N_1) = \max \{s_1,
\dots, s_{N_1}\} = \bfa{1}_M(\x_1)$. Therefore, the statement is true for
$s = 2$ and $r = 0$. Assume that either $s > 2$ or $r > 0$. It can
be seen that $v^\brak{s}(0) = v^\brak{s-1}(N_{s-1})$. Thus, we can
assume that $r > 0$ (and $s \ge 2$).

Let $N = M/(\x_1, \dots, \x_{s-1}, x_{s,1}, x_{s,2}, \dots,
x_{s,r-1})M$. Observe first that $N_\m \not= 0$ for some $\m$ if and only if $m_1$ equals the first coordinate of a
generating degree of $M$ (i.e., $m_1 = d_{i,1}(M)$ for some $i$).
This implies that
\begin{align}
\a{1}(0 :_N x_{s,r}) \le \a{1}(N) = d^\brak{1}(M). \label{0piece}
\end{align}

Consider an arbitrary $\n = (n_1, \dots, n_k) \in \Z^k$ with $n_1 >
\max\{v^\brak{s}(r-1), d^\brak{1}(M)\}$. By Lemma \ref{homology
exact sequence}, we have
\begin{align}
& H^\crak{s}_i(r-1;M)_{\n+i\e_1} \rightarrow
H^\crak{s}_i(r;M)_{\n+i\e_1} \rightarrow H^\crak{s}_{i-1}(r-1;M)_{\n+i\e_1 - \e_s} \label{seqr2} \\
& H^\crak{s}_1(r-1;M)_{\n+\e_1} \rightarrow
H^\crak{s}_1(r;M)_{\n+\e_1} \rightarrow
\tilde{H}^\crak{s}_0(r-1;M)_{\n+\e_1-\e_s} \rightarrow 0.
\label{seqr3}
\end{align}

Since $n_1 > \max \{v^\brak{s}(r-1), d^\brak{1}(M)\}$, it follows
from (\ref{0piece}) and (\ref{seqr3}) that
$$H^\crak{s}_1(r;M)_{\n+\e_1} = 0.$$ Since $n_1 > v^\brak{s}(r-1)$, we have $H^\crak{s}_i(r-1;M)_{\n+i\e_1} =
H^\crak{s}_{i-1}(r-1;M)_{\n+i\e_1-\e_s} = 0$. Thus, (\ref{seqr2})
implies that $H^\crak{s}_i(r;M)_{\n+i\e_1} = 0$ for all $i \ge 2$.
Hence,
\begin{align}
v^\brak{s}(r) \le \max \{v^\brak{s}(r-1), d^\brak{1}(M)\}.
\label{1st inequality}
\end{align}

Observe now that if $v^\brak{s}(r-1) \le d^\brak{1}(M)$ then
(\ref{1st inequality}) implies that
$$\max \{v^\brak{s}(r), d^\brak{1}(M)\} = d^\brak{1}(M) = \max \{v^\brak{s}(r-1), d^\brak{1}(M)\}.$$
On the other hand, if $v^\brak{s}(r-1) > d^\brak{1}(M)$ then there
exist an integer $i \ge 1$ and $\n \in \Z^k$
with $n_1 = v^\brak{s}(r-1)$ such that
$H^\crak{s}_i(r-1;M)_{\n+i\e_1} \not= 0$. In this case, if
$H^\crak{s}_i(r;M)_{\m+i\e_1} = 0$ for any $\m \in \Z^k$ with $m_1 \ge v^\brak{s}(r-1)$, then by Lemma
\ref{homology exact sequence} we have
$$0 \rightarrow \bigoplus_{m_1 = n_1} H^\crak{s}_i(r-1;M)_{\m+i\e_1-\e_s} \stackrel{x_{s,r}}{\rightarrow} \bigoplus_{m_1 = n_1} H^\crak{s}_i(r-1;M)_{\m+i\e_1} \rightarrow 0.$$
This implies that
$$\bigoplus_{m_1 = n_1} H^\crak{s}_i(r-1;M)_{\m+i\e_1} = x_{s,r} \Big(\bigoplus_{m_1 = n_1} H^\crak{s}_i(r-1;M)_{\m+i\e_1} \Big),$$
which is a contradiction by Nakayama's lemma. Therefore, there
exists, in this case, $\m_0 = (m_{01}, \dots, m_{0k}) \in \Z^k$ with
$m_{01} = v^\brak{s}(r-1)$ such that $H^\crak{s}_i(r;M)_{\m_0+i\e_1}
\not= 0$. That is, $v^\brak{s}(r) \ge v^\brak{s}(r-1) >
d^\brak{1}(M)$. We, hence, have $\max \{v^\brak{s}(r),
d^\brak{1}(M)\} = \max \{v^\brak{s}(r-1), d^\brak{1}(M)\}$. The
conclusion now follows by induction.
\end{proof}

\begin{example} \label{RunExample3}
Consider $S = \kk[x,y,z]$ and $M = S/(x^2y,xy^2,xyz, y^3,
y^2z,yz^2)$ as in Example \ref{RunExample1}. Notice that $\dim M =
2$. It can be seen that $\bfa{1}_M(\langle x,z \rangle)=1$, $\bfa{2}_M(\langle x,z \rangle)=2$, and
$\bfa{3}_M(\langle x,z \rangle)=1$ while
$d^\brak{1}(M)=d^\brak{2}(M)=d^\brak{3}(M)=0$. Theorem \ref{reg:
filter-regular} implies that $\Reg(M)=(1,2,1)$.
\end{example}

The following lemma (see \cite[Lemma 2.3]{RedExp}) exhibits the behavior of modding out by a filter-regular element with respect a given coordinate.

\begin{lemma} \label{1filter-regular}
Let $g \in S_{\e_l}$ be an $M$-filter-regular element with respect to the
$l$-th coordinate and let $(\x_l')$ be the image of the ideal
$(\x_l)$ in $S/gS$. Then for all $i \ge 0$,
\begin{align*}
\a{l}\big( H^{i+1}_{(\x_l)}(M) \big) + 1 & \le \a{l} \big(
H^i_{(\x_l')}(M/
gM) \big) \\
& \le \max \big\{ \a{l}\big(H^i_{(\x_l)}(M)\big),
\a{l}\big(H^{i+1}_{(\x_l)}(M)\big) + 1\big\}.
\end{align*}
\end{lemma}

\begin{proof} It follows from the proof of Lemma \ref{primes avoidance} that $(0 :_M g) \subseteq \bigcup_{n=0}^\infty [0 :_M (\x_l)^n].$ That is, $(0 :_M g)$ is annihilated by some power of $(\x_l)$. This implies that
\begin{align}
H^i_{(\x_l)}(0:_M g) = 0 \ \forall i \ge 1. \label{1filter-eq1}
\end{align}

Consider the exact sequence
$$0 \rightarrow (0 :_M g) \rightarrow M \rightarrow M/(0 :_M g) \rightarrow 0.$$
(\ref{1filter-eq1}) implies that $H^i_{(\x_l)}(M) =
H^i_{(\x_l)}(M/(0:_M g))$ for all $i \ge 1$. Now, consider the long
exact sequence of cohomology groups associated to the exact sequence
$$0 \rightarrow M/(0 :_M g)(-\e_l) \stackrel{g}{\rightarrow} M \rightarrow
M/gM \rightarrow 0$$ we get
\begin{align}
H^i_{(\x_l)}(M)_\n \rightarrow H^i_{(\x_l')}(M/gM)_\n \rightarrow
H^{i+1}_{(\x_l)}(M)_{\n-\e_l} \stackrel{g}{\rightarrow}
H^{i+1}_{(\x_l)}(M)_\n \label{1filter-eq2}
\end{align}
for any $i \ge 0$ and $\n \in \Z^k$.

By the definition of the function $\a{l}$, there exists $\m \in \Z^k$ with $m_l =
\a{l}\big(H^{i+1}_{(\x_l)}(M)\big)+1$ such that
$H^{i+1}_{(\x_l)}(M)_\m = 0$ and $H^{i+1}_{(\x_l)}(M)_{\m-\e_l}
\not= 0$. (\ref{1filter-eq2}) then implies that
$H^i_{(\x_l')}(M/gM)_\m \not= 0$. Thus,
$$\a{l}\big(H^i_{(\x_l')}(M/gM)\big) \ge \a{l}\big(H^{i+1}_{(\x_l)}(M)\big)+1.$$ On the other hand, for any
$\n \in \Z^k$ such that
$$n_l > \max \big\{\a{l}\big(H^i_{(\x_l)}(M)\big), \a{l}\big(H^{i+1}_{(\x_l)}(M)\big)+1\big\},$$ we have
$H^i_{(\x_l)}(M)_\n = H^{i+1}_{(\x_l)}(M)_{\n-\e_l} = 0$.
(\ref{1filter-eq2}) now implies that
$$\a{l} \big( H^i_{(\x_l')}(M/gM) \big) \le \max \big\{ \a{l}\big(H^i_{(\x_l)}(M)\big), \a{l}\big(H^{i+1}_{(\x_l)}(M)\big) + 1\big\}.$$
The lemma is proved.
\end{proof}

The next theorem lays the groundwork for the rest of the paper,
allowing us to prove our main theorems in Section~\ref{s.asympreg}.

\begin{theorem} \label{reg: cohomology and end}
Let $S$ be a standard $\N^k$-graded polynomial ring over a field $\kk$ and
let $M$ be a finitely generated $\Z^k$-graded $S$-module. Suppose that $\x_l$ is an
$M$-filter-regular sequence with respect to the $l$-th coordinate. Then,
\begin{enumerate}
\item $\bfa{l}_M(\x_l) = \max \Big\{ \a{l}\big(H^i_{(\x_l)}(M)\big) + i ~\big|~ 0 \le i \le N_l - 1 \Big\}.$
\item $\bfa{l}_M(\x_l) = \max \Big\{ \a{l} \Big( \dfrac{(x_{l,1}, \dots, x_{l,i})M :_M (\x_l)}{(x_{l,1}, \dots, x_{l,i})M} \Big) ~\Big|~ 0 \le i \le N_l-1 \Big\}.$
\item $\reg_l(M) = \max \Big\{ \a{l} \Big( \dfrac{(x_{l,1}, \dots, x_{l,i})M :_M (\x_l)}{(x_{l,1}, \dots, x_{l,i})M} \Big) ~\Big|~ 0 \le i \le N_l \Big\}.$
\end{enumerate}
\end{theorem}

\begin{proof} We prove (1) using induction on $N_l$. Since $x_{l,1}$ is an
$M$-filter-regular element with respect to the $l$-th coordinate, as in the
proof of Lemma \ref{primes avoidance}, we have $(0:_M x_{l,1})
\subseteq \bigcup_{n=0}^\infty [0 :_M (\x_l)^n] = H^0_{(\x_l)}(M)$.
Thus, $\bfa{l}_M(x_{l,1}) \le \a{l}\big(H^0_{(\x_l)}(M)\big)$. On
the other hand, there exists $\n \in \Z^k$ with
$n_l = \a{l}\big(H^0_{(\x_l)}(M)\big)$ such that $H^0_{(\x_l)}(M)_\n
\not= 0$. Furthermore, $x_{l,1} H^0_{(\x_l)}(M)_\n \subseteq
H^0_{(\x_l)}(M)_{\n+\e_l} = 0$, which implies that
$H^0_{(\x_l)}(M)_\n \subseteq (0 :_M x_{l,1})_\n$. Therefore,
$\bfa{l}_M(x_{l,1}) \ge \a{l}\big(H^0_{(\x_l)}(M)\big)$. Hence,
$\bfa{l}_M(x_{l,1}) = \a{l}\big(H^0_{(\x_l)}(M)\big)$ and the
statement is true for $N_l = 1$.

Assume that $N_l > 1$. Let $\x_l'$ denote the sequence of images of
$x_{l,2}, \dots, x_{l,N_l}$ in $S/x_{l,1}S$ and let $N =
M/x_{l,1}M$. Then $\x_l'$ gives an $N$-filter-regular sequence with respect to
the $l$-th coordinate. By induction and using Lemma
\ref{1filter-regular}, we have
\begin{align*}
\lefteqn{\max\{\a{l}\big(H^i_{(\x_l)}(M)\big)+i ~|~ i = 1, \dots,
N_l-1\}} \\
& \le \max \{ \a{l}\big(H^{i-1}_{(\x_l')}(N)\big)+i-1 ~|~ i = 1,
\dots, N_l-1\} = \bfa{l}_N(\x_l') \\
& \le \max \{\a{l}\big(H^j_{(\x_l)}(M)\big)+j ~|~ j = 0, \dots,
N_l-1\}.
\end{align*}
This, together with the fact that $\bfa{l}_M(\x_l) = \max
\{\bfa{l}_M(x_{l,1}), \bfa{l}_N(\x_l')\}$, implies that
$$\bfa{l}_M(\x_l) = \max \big\{\a{l}\big(H^i_{(\x_l)}(M)\big)+i ~|~ i = 0, \dots, N_l-1\big\}.$$
(1) is proved.

To prove (2), we first successively apply Lemma
\ref{1filter-regular} to get
\begin{align*}
\a{l}\big(H^i_{(\x_l)}(M)\big)+i & \le
\a{l}\big(H^0_{(\x_l)S_i}(M/(x_{l,1}, \dots, x_{l,i-1})M)\big) \\
& \le \max \big\{\a{l}\big(H^j_{(\x_l)}(M)\big)+j ~|~ j = 0, \dots,
i \big\}
\end{align*}
for any $i \ge 0$, where $S_i = S/(x_{l,1}, \dots, x_{l,i-1})S$. It
follows that for any $t \le N_l$,
\begin{align}
\max_{0 \le i \le t} \big\{ \a{l}\big(H^i_{(\x_l)}(M)\big)+i \big\}
= \max_{0 \le i \le t} \big\{ \a{l}\big( H^0_{(\x_l)S_i}(M/(x_{l,1},
\dots, x_{l,i-1})M)\big) \big\}. \label{regcoend-eq1}
\end{align}

Fix $i$ and consider an arbitrary $\m \in \Z^k$ such that
$$m_l = \a{l}\big(H^0_{(\x_l)S_i}(M/(x_{l,1}, \dots,
x_{l,i-1})M)\big).$$ Then, $(\x_l)
H^0_{(\x_l)S_i}(M/(x_{l,1},\dots,x_{l,i-1})M)_\m \subseteq
H^0_{(\x_l)S_i}(M/(x_{l,1},\dots,x_{l,i-1})M)_{\m+\e_l} = 0$.
Moreover, $H^0_{(\x_l)S_i}(M/(x_{l,1}, \dots, x_{l,i-1})M) =
\bigcup_{n=0}^\infty \Big(\dfrac{(x_{l,1}, \dots, x_{l,i-1})M :_M
(\x_l)}{(x_{l,1}, \dots, x_{l,i-1})M}\Big)$. Thus, we have
\begin{align*}
H^0_{(\x_l)S_i}(M/(x_{l,1}, \dots, x_{l,i-1})M)_\m & \subseteq \dfrac{(x_{l,1}, \dots, x_{l,i-1})M :_M (\x_l)}{(x_{l,1}, \dots, x_{l,i-1})M} \\
& \subseteq H^0_{(\x_l)S_i}(M/(x_{l,1}, \dots, x_{l,i-1})M).
\end{align*}
This implies that
\begin{align}
\a{l}\Big(\dfrac{(x_{l,1}, \dots, x_{l,i-1})M :_M (\x_l)}{(x_{l,1},
\dots, x_{l,i-1})M}\Big) = m_l. \label{regcoend-eq2}
\end{align}
(2) now follows from part (1), (\ref{regcoend-eq1}) and
(\ref{regcoend-eq2}).

In view of Theorem \ref{reg: filter-regular} and part (2), to prove
(3) we only need to show that
$$\a{l}(M/(\x_l)M) = d^\brak{l}(M).$$
This is indeed true since $\bigoplus_{n_l=a} [M/(\x_l)M]_{(n_1,
\dots, n_k)} \not= 0$ if and only if $a$ equals to the $l$-th
coordinate of a generating degree of $M$ (i.e. $a = d_{i,l}(M)$ for
some $i$). (3) is proved.
\end{proof}


\section{Asymptotic behavior of resolution regularity} \label{s.asympreg}

In this section, we prove our main theorems. Having Theorem
\ref{reg: cohomology and end}, our arguments now are direct
generalization to the $\Z^k$-graded situation of techniques used in
\cite{TW}.

Our first main result establishes the asymptotic linearity of $\Reg(I^nM)$.

\begin{theorem} \label{asymptotic linearity}
Let $S$ be a standard $\N^k$-graded polynomial ring over a field $\kk$ and
let $M$ be a finitely generated $\Z^k$-graded $S$-module. Let $I$ be
a multigraded homogeneous ideal in $S$ minimally generated in
degrees $\d_1(I ), \dots, \d_v(I)$. Then, $\Reg(I^nM)$ is
asymptotically a linear function with slope vector at most $\big(
d^\brak{1}(I), \dots, d^\brak{k}(I) \big)$ componentwise.
\end{theorem}

\begin{proof} Observe first that if for each $1 \le l \le k$, $\reg_l(I^nM)$ is asymptotically a linear function $a_ln + b_l$, then $\Reg(I^nM)$ is asymptotically the linear function $(a_1, \dots, a_k)n + (b_1, \dots, b_k)$. Thus, it suffices to show that for each $1 \le l \le k$, $\reg_l(M)$ is asymptotically indeed a linear function with slope at most $d^\brak{l}(I)$.

Assume that $\{F_1, \dots, F_v\}$ is a minimal set of generators for $I$, where $\deg
F_i = \d_i(I)$. Let $\M = \bigoplus_{n \ge 0}I^nM$ be the Rees
module of $M$ with respect to $I$. Then $\M$ is a finitely generated
$\Z^{k+1}$-graded module over the Rees algebra of $I$, $\R =
\bigoplus_{n \ge 0}I^nt^n$. Observe further that there is a natural
surjection of $\Z^{k+1}$-graded algebras $R=\kk[\X_1, \dots, \X _k,
Y_1, \dots, Y_v] \rightarrow \R$ given by $\X_i \mapsto \x_i$ and
$Y_s \mapsto F_st$, where $\deg X_{i,j} = (\e_i, 0) \in \Z^{k+1}$
and $\deg Y_s = (\d_s(I), 1) \in \Z^{k+1}$ for all $1 \le i \le k$,
$1 \le j \le N_i$ and $1 \le s \le v$.

Under this surjection, $\M$ can be viewed as a finitely generated
$\Z^{k+1}$-graded $R$-module. Moreover, $\M_n = I^nM$.
The conclusion now follows from Theorem \ref{linear: regularity} below.
\end{proof}

Let $\X_i$ represent the set of indeterminates $\{X_{i,1}, \dots,
X_{i,N_i}\}$ for $i = 1, \dots, k$. Let $\d_1, \dots, \d_w \in
\Z^k$, where $\d_i = (d_{i,1}, \dots, d_{i,k})$, be non-negative
vectors (i.e. $d_{i,j} \ge 0$ for all $i$ and $j$). For each $1 \le
l \le k$, let
$$d^\brak{l} = \max \{d_{i,l} ~|~ 1 \le i \le w\} \ge 0.$$

\begin{theorem} \label{linear: regularity}
Let $R = \kk[\X_1, \dots, \X_k, Y_1, \dots, Y_w]$ be an
$\N^{k+1}$-graded polynomial ring where $\deg X_{i,j} = (\e_i, 0)
\in \Z^{k+1}$ for all $i=1, \dots, k$ and $j=1, \dots, N_i$, and
$\deg Y_s = (\d_s, 1) \in \Z^{k+1}$ for all $s = 1, \dots, w$. Let
$\M = \bigoplus_{\k \in \Z^{k+1}} \M_\k$ be a finitely generated
$\Z^{k+1}$-graded $R$-module. For each $n \in \Z$, let
$$\M_n = \bigoplus_{\m \in \Z^k} \M_{(\m,n)}.$$
Then for each $1 \le l \le k$, $\reg_l(\M_n)$ is asymptotically a
linear function with slope at most $d^\brak{l}$.
\end{theorem}

\begin{proof} Considering the flat base extension $\kk \rightarrow \kk' =
 \kk(U)$ where $U = (u_{ij})_{1 \le i,j \le N_l}$ is a square matrix of indeterminates. By Lemma \ref{flat extension} and a change of variables, we may assume that $\x_l$ is
an $\M_n$-filter-regular sequence with respect to the $l$-th coordinate for
any $n \in \Z$. By Theorem \ref{reg: cohomology and end}, we have
\begin{align}
\reg_l(\M_n) = \max \Big\{ \a{l} \Big( \dfrac{(x_{l,1}, \dots,
x_{l,i})\M_n :_{\M_n} (\x_l)}{(x_{l,1}, \dots, x_{l,i})\M_n} \Big)
~\Big|~ 0 \le i \le N_l \Big\}. \label{linear: regularity-eq1}
\end{align}

It is easy to see that
$$\dfrac{(x_{l,1}, \dots, x_{l,i})\M_n :_{\M_n} (\x_l)}{(x_{l,1}, \dots, x_{l,i})\M_n} = \Big[\dfrac{(x_{l,1}, \dots, x_{l,i})\M :_{\M} (\x_l)}{(x_{ l,1}, \dots, x_{l,i})\M} \Big]_n.$$ Moreover, we can
consider $\dfrac{(x_{l,1}, \dots, x_{l,i})\M :_{\M}
(\x_l)}{(x_{l,1}, \dots, x_{l,i})\M}$ as a $\Z^k$-graded module over the polynomial ring
$\kk[\X_1, \dots, \widehat{\X_l}, \dots, \X_k, Y_1, \dots, Y_w],$
where $\widehat{\X_l}$ indicates that $\X_l$ is removed. It remains to show that
$\a{l}\Big(\Big[\dfrac{(x_{l,1}, \dots, x_{l,i})\M :_{\M} (\x_l)}{(x_{l,1}, \dots, x_{l,i})\M}\Big]_n\Big)$ is asymptotically a linear function with slope
$\le d^\brak{l}$ for each $i = 0, \dots, N_l$. Indeed, this follows from a more general statement of our next result, Theorem \ref{linear: a-invariant}.
\end{proof}

\begin{theorem} \label{linear: a-invariant}
Let $0 \le t \le k$ and let $1 \le l_1 < \dots < l_t \le k$. Let
$$R = \kk[\X_{l_1}, \dots, \X_{l_t}, Y_1, \dots, Y_w]$$
be a $\N^{k+1}$-graded polynomial ring where $\deg X_{l_i,j} =
(\e_{l_i}, 0) \in \Z^{k+1}$ for all $i=1, \dots, t$ and $j=1, \dots,
N_{l_i}$, and $\deg Y_s = (\d_s, 1) \in \Z^{k+1}$ for all $s = 1,
\dots , w$. Let $\M = \bigoplus_{\k \in \Z^{k+1}} \M_\k$ be a
finitely generated $\Z^{k+1}$-graded $R$-module. For each $n \in
\Z$, let $$\M_n = \bigoplus_{\m \in \Z^k} \M_{(\m,n)}.$$ Then for
each $1 \le l \le k$, $\a{l}(\M_n)$ is asymptotically a linear
function with slope at most $d^\brak{l}$.
\end{theorem}

\begin{proof} We use induction on $w$. Suppose $w = 0$. Let $h = d^\brak{k+1}(\M)$ be the maximum of the $(k+1)$-st coordinate of minimal
generating degrees of $\M$. Then, it is easy to see that for $n \gg 0$, specifically for $ n > h$, we must have $\M_n = 0$. Thus, $\a{l}(\M_n)$ is asymptotically the zero function.

Assume now that $w \ge 1$. Consider the exact sequence of
$\Z^{k+1}$-graded $S$-modules
$$0 \rightarrow [0 :_\M Y_w]_{(\m,n)} \rightarrow \M_{(\m,n)} \stackrel{Y_w}{\rightarrow} \M_{(\m+\d_w, n+1)} \rightarrow [\M/Y_w\M]_{(\m+\d_w,n+1)} \rightarrow 0.$$
Since $0 :_\M Y_w$ and $\M/Y_w\M$ can be viewed as $\Z^{k+1}$-graded
modules over the ring $\kk[\X_{l_1}, \dots, \X_{l_t}, Y_1, \dots, Y_{w-1}]$, by induction $\a{l}([0 :_\M Y_w]_n)$ and $\a{l}([\M/Y_w\M]_n)$ are asymptotically linear functions with
slopes at most $d^\brak{l}$. As a consequence, we have
\begin{align}
\a{l}([0 :_\M Y_w]_n) + d^\brak{l} & \ge \a{l}([0 :_\M Y_w]_{n+1}),
\text{and} \label{ineq1} \\
\a{l}([\M/Y_w\M]_n) + d^\brak{l} & \ge \a{l}([\M/Y_w\M]_{n+1})
\label{ineq2}
\end{align}
for all $n \gg 0$.

It is clear that $\a{l}(\M_n) \ge \a{l}([0 :_\M
Y_w]_n)$ for all $n \ge 0$. If $\a{l}(\M_n) = \a{l}([0 :_\M Y_w]_n)$
for $n \gg 0$ then the statement follows by induction. It remains to
consider the case that there exists an infinite sequence of integers
$n$ for which $\a{l}(\M_n) > \a{l}([0 :_\M Y_w]_n)$.

Let $m$ be such an integer. Making use of the inequalities
(\ref{ineq1}) and (\ref{ineq2}), we obtain
$$\a{l}(\M_{m+1}) = \max \{ \a{l}(\M_m) + d^\brak{l}, \a{l}([\M/Y_w\M]_{m+1}) \}$$ (since we can take $m$ to be bigger than the last
coordinate of all generators in a minimal system of generators for
$\M$). By taking $m$ large enough, we further have
$$\a{l}(\M_m) + d^\brak{l} \ge \a{l}([\M/Y_w\M]_m) + d^\brak{l} \ge \a{l}([\M/Y_w\M]_{m+1}).$$
Thus, $\a{l}(\M_{m+1}) = \a{l}(\M_m) + d^\brak{l} > \a{l}([0 :_\M
Y_w]_m) + d^\brak{l} \ge \a{l}([0 :_\M Y_w]_{m+1})$.

We have shown that for $m \gg 0$, if $\a{l}(\M_m) > \a{l}([0 :_\M Y_w]_m)$ then $\a{l}(\M_{m+1}) > \a{l}([0 :_\M Y_w]_{m+1})$.
By repeating this argument we can show that $\a{l}(\M_n) >
\a{l}([0:_\M Y_w]_n)$ for all $n \gg 0$. This, by a similar argument
as above, implies that $\a{l}(\M_{n+1}) = \a{l}(\M_n) + d^\brak{l}$
for all $n \gg 0$. Hence, $\a{l}(\M_n)$ is also asymptotically a linear
polynomial with slope $d^\brak{l}$ in this case. The assertion, and hence, Theorem \ref{linear: regularity} and subsequently Theorem \ref{asymptotic linearity}, is proved.
\end{proof}

The rest of this section is to explicitly describe the slope vector
of the asymptotically linear function $\Reg(I^nM)$.

Recall from Definition~\ref{def.reduction} that
a multigraded homogeneous ideal $J \subset I$ is an {\it $M$-reduction} of $I$ if $I^nM = JI^{n-1}M$ for all $n \gg 0$.
For each $l = 1, \dots, k$,
$$\rho^\brak{l}_M(I) = \min \{ d^\brak{l}(J) ~\big|~ J \mbox{ is an } M-\mbox{reduction of } I\},$$ and
$$\beg^\brak{l}(M) = \min \{ n_l ~\big|~ \exists \n \in \Z^k \mbox{ such that } M_\n \not= 0\}.$$
It can further be seen that if $I^{n_0}M = JI^{n_0-1}M$ then $I^nM = JI^{n-1}M$ for
all $n \ge n_0$ and $J$ is an $M$-reduction of $I$.

\begin{lemma} \label{Degree Lower Bound}
For any $n \ge 0$, we have $$d^\brak{l}(I^nM) \ge n
\rho^\brak{l}_M(I) + \beg^\brak{l}(M).$$
\end{lemma}

\begin{proof} Write $I$ as the sum of two ideals, $I = J+K$, in which $J$ is generated by multigraded homogeneous elements $f$ of $I$ such that $\deg_l(f) < \rho^\brak{l}_M(I)$ and $K$ is generated by multigraded
homogeneous elements $g$ of $I$ such that
$\deg_l(g)\ge\rho^\brak{l}_M(I)$. Then,
$$I^nM = (J+K)I^{n-1}M = JI^{n-1}M + K(J+K)^{n-1}M = JI^{n-1}M + K^nM.$$
It can be seen that $K^nM$ is generated by multigraded homogeneous elements
$h$ for which $\deg_l(h) \ge n \rho^\brak{l}_M(I) +
\beg^\brak{l}(M)$.

If for some $n_0 \ge 0$, $d^\brak{l}(I^{n_0}M) < n_0
\rho^\brak{l}_M(I) + \beg^\brak{l}(M)$ then $K^{n_0}M = 0$ and
$I^{n_0}M = JI^{n_0-1}M$. This implies that $I^nM = JI^{n-1}M$ for
any $n \ge n_0$ and $J$ is an $M$-reduction of $I$. However, by the
construction of $J$, we have $d^\brak{l}(J) < \rho^\brak{l}_M(I)$,
which is a contradiction to the definition of $\rho ^\brak{l}_M(I)$.
The statement is proved.
\end{proof}

The next theorem is a refinement of Theorem~\ref{asymptotic
linearity} and gives an explicit description for the slope vector of
the asymptotically linear function $\Reg(I^nM)$.

\begin{theorem} \label{coefficient}
Let $S$ be a standard $\N^k$-graded polynomial ring over a field $\kk$ and
let $M$ be a finitely generated $\Z^k$-graded $S$-module. Let $I$ be
a multigraded homogeneous ideal in $S$. Then, there exists a vector $\b \ge
(\beg^\brak{1}(M), \dots, \beg^\brak{ k}(M))$ such that for $n \gg
0$,
$$\Reg(I^nM) = n \big(\rho^\brak{1}_M(I), \dots, \rho^\brak{k}_M(I)\big)
+ \b.$$
\end{theorem}

\begin{proof} Suppose $1 \le l \le k$. Let a multigraded homogeneous ideal
$J$ be an $M$-reduction of $I$ such that $d^\brak{l}(J) =
\rho^\brak{l}_M (I)$. Let $\M$ be the Rees
module of $M$ with respect to $I$ and let $\R$ be the Rees algebra of $I$. As before, $\M$ is a finitely
generated $\Z^{k+1}$-graded module over the Rees algebra of $I$, $\R
= \bigoplus_{n \ge 0}I^nt^n$. Let $S[Jt]$ be the Rees algebra of
$J$. Then, since $J$ is an $M$-reduction of $I$, $S[Jt] \rightarrow
\R$ is a finite map. It follows that $\M$ is a finitely generated
$\Z^{k+1}$-graded module over $S[Jt]$.

Assume that $J$ is
minimally generated by $G_1, \dots, G_u$ of degrees $\d_1(J), \dots,
\d_u(J)$. Then, there is a natural surjection $R = \kk[\X_1, \dots,
\X_k, Y_1, \dots, Y_u] \rightarrow S[Jt]$ given by $\X_i \mapsto
\x_i$ and $Y_s \mapsto G_st$, where $\deg X_{i,j} = (\e_i,0) \in
\Z^{k+1}$ and $\deg Y_s = (\d_s(J),1) \in \Z^{k+ 1}$ for all $1 \le
i \le k$, $1 \le j \le N_i$ and $1 \le s \le u$. Under this
surjection, $\M$ can be viewed as a finitely generated
$\Z^{k+1}$-graded module over $R$. By Theorem~\ref{linear:
regularity}, for each $l = 1, \dots, k$, $\reg_l(I^nM) =
\reg_l(\M_n)$ is asymptotically a linear function with slope at most
$d^\brak{l}(J)$. Suppose, $\reg_l(I^nM) = a_l n + b_l$ for $n \gg
0$, in which $a_l \le d^\brak{l}(J)$.

By Lemma \ref{Degree Lower Bound}, we  have
\begin{align}
\reg_l(I^nM) \ge d^\brak{l}(I^nM) \ge n \rho^\brak{l}_M(I) +
\beg^\brak{l}(M) \ \forall \n \ge 0. \label{eq525}
\end{align}
This implies that $a_l \ge \rho^\brak{l}_M(I)$. Together with the
fact that $a_l \le d^\brak{l}(J) = \rho^\brak{l}_M(I)$, we now have
$a_l = \rho^\brak{l}_M(I)$. (\ref{eq525}) also implies that $b_l
\ge\beg^\brak{l}(M)$. The theorem is proved.
\end{proof}

\begin{example} Let $S=\mathbb{Q}[a,b,c]$ be a standard $\N^3$-graded polynomial ring, where $\deg a = \e_1, \deg b = \e_2$ and $\deg c = \e_3$. Consider a homogeneous ideal $I = (a,b,c)$ and a $\Z^3$-graded module $M = S/(a^3,a^2b,a^2c,abc)$.
Observe that $J = (b,c)$ is a minimal $M$-reduction of $I$. Thus,
$(\rho^\brak{1}_M(I), \rho^\brak{2}_M(I), \rho^\brak{3}_M(I)) =
(0,1,1) \not= (1,1,1) = (d^\brak{1}(I), d^\brak{2}(I), d^\brak{3}(I))$. Moreover, $(\beg^\brak{1}(M), \beg^\brak{2}(M),
\beg^\brak{3}(M)) = \0$. Theorem~\ref{coefficient} says that there
exists $\b \ge \0$ so that for $n \gg 0$ we have
$$\Reg(I^nM) = n(0,1,1) + \b.$$

Using CoCoA \cite{cocoa}, we can also calculate $\Reg(IM) =
(1,1,1)$, $\Reg(I^2M) = (2,2,2) \not= 2(0,1,1) + (1,0,0)$, while $\Reg(I^3M) = 3(0,1,1) +
(1,0,0)$, $\Reg(I^4M) = 4(0,1,1) + (1,0,0)$, and $\Reg(I^5M) =
5(0,1,1) + (1,0,0)$. It seems that in this example, $\Reg(I^nM) = n(0,0,1) + (1,0,0)$ for $n \ge 3$.
\end{example}

\begin{remark} It would be interesting to find the vector $\b$ and a lower bound $n_0$ so that
in Theorem \ref{coefficient},
$$\Reg(I^nM) = n(\rho^\brak{1}_M(I), \dots, \rho^\brak{k}_M(I)) + \b$$
for all $n \ge n_0$. Even in the $\Z$-graded case, this is still
open.
\end{remark}


\section{Appendix} \label{s.infinite}

During the preparation of this paper, Bernd Ulrich asked us whether
multigraded regularity can be defined even when a {\it finite}
minimal free resolution as in Definition \ref{def.reg} may not exist.
Ulrich's question has helped us better understand
resolution regularity. We would like to thank him for his question. We shall also propose an alternative definition for multigraded regularity even when finite minimal free resolutions may not exist.

Let $A$ be a Noetherian ring (commutative with identity). Let $T$ be
a standard $\Z^k$-graded algebra over $A$; that is, $T_\0 = A$ and
$T$ is generated over $T_\0$ by elements of $\bigoplus_{l=1}^k
T_{\e_l}$. Let $N = \bigoplus_{\n \in \Z^k}N_\n$ be a finitely
generated $\Z^k$-graded module over $T$. It is well known that, in
general, a finite minimal free resolution of $N$ over $T$ may not
exist. Thus, Definition \ref{def.reg} will not work in this case. We
shall introduce an alternative approach to define $\nReg(N )$. This
approach results from our Theorem \ref{reg: cohomology and end}.
Our new notion, $\nReg(N)$, still bounds (but does not capture) the supremum of coordinates (with suitable shifts) in a minimal free resolution of $N$.

As before, for each $l = 1, \dots, k$, assume that $T_{\e_l}$ is
generated as a module over $T_\0$ by $\x_l = \{x_{l,1}, \dots,
x_{l,N_l}\}$ . Write $T = \bigoplus_{m \in N}T_m$, where $T_m =
\bigoplus_{\n \in \N^k, n_l = m} T_\n$. Then,
$T$ can be viewed as a standard $\N$-graded algebra over $T_0 =
\bigoplus_{\n \in \N^k, n_l = 0} T_\n$. Under
this new grading of $T$, $(\x_l)$ is the homogeneous irrelevant
ideal. We can also view $N$ as a $\Z$-graded $T$-module with
the grading given by $N = \bigoplus_{m \in \Z} N_m$, where $N_m =
\bigoplus_{\n \in \Z^k, n_l = m} N_\n$. This
induces a natural $\Z$-graded structure on local cohomology modules
$ H^i_{(\x_l)}(N)$ for all $i \ge 0$. Under this $\Z$-graded structure (depend on $l$), the usual $a$-invariants are well defined and have been much studied (cf. \cite{Sh, trung}). We shall use the superscript $\bullet^\brak{l}$ to indicate that these invariants depend on $l$. That is, we define
$$a^\brak{l}_i(N) = \left\{ \begin{array}{lll} -\infty & \mbox{if} & H^i_{(\x_l)}(N) = 0 \\
\max\{ m ~|~ H^i_{(\x_l)}(N)_m \not= 0\} & \mbox{if} &
H^i_{(\x_l)}(N) \not= 0.
\end{array} \right.$$

\begin{definition} \label{def.alt}
For each $l = 1, \dots, k$, let
$$\nreg_l(N) = \max \{ a^\brak{l}_i(N) + i ~|~ i \ge 0\}.$$
We define the multigraded (resolution) {\it regularity} of $N$ to be the vector
$$\nReg(N) = (\nreg_1(N), \dots, \nreg_k(N)) \in \Z^k.$$
\end{definition}

Observe that since $H^i_{(\x_l)}(N) = 0$ for all $i \gg 0$,
$\nReg(N)$ is well defined.

\begin{remark}
In our current setting, $N$-filter regular
sequences with respect to the $l$-th coordinate, and invariants
$\bfa{l}_N(\bullet)$ and $\a{l}(\bullet)$ can be defined as before.
By applying a similar line of arguments as that of Theorem~\ref{reg: cohomology and end}.(1), we have $$\nreg_l(N) = \max \{ \bfa{l}_N(\x_l), \a{l}\big(H^{N_l}_{(\x_l)}(N)\big) + N_l \}.$$
This, together with a similar line of arguments as that of Theorem~\ref{reg: cohomology and end}.(2), implies that
\begin{align}
 \nreg_l(N) & = \max \{ \bfa{l}_N(\x_l), d^\brak{l}(N) \} \label{app: reg and a-inv} \\
 & = \max \Big\{ \a{l} \Big( \dfrac{(x_{l,1}, \dots,
x_{l,i})N :_N ( \x_l)}{(x_{l,1}, \dots, x_{l,i})N} \Big) ~\Big|~ 0 \le i \le N_l \Big\}. \label{cohend}
\end{align}

Let $S \rightarrow T$ be the natural surjection from a standard $\N^k$-graded polynomial ring over $A$ to $T$. Then $N$ is a $\Z^k$-graded $S$-module.
Let $\GG: \dots \rightarrow G_p \rightarrow \dots \rightarrow G_0 \rightarrow N \rightarrow 0$
be a (infinite) minimal free resolution of $N$ over $S$, where $G_i = \bigoplus_j S(-c_{ij}^\brak{1}, \dots, -c_{ij}^\brak{k})$ for all $i \ge 0$. For each $1 \le l \le k$ and $i \ge 0$, set $c_i^\brak{l} = \max_j \{c_{ij}^\brak{l}\}$, and let
$$r_l = \sup_i \{ c_i^\brak{l} - i \}.$$
Let $\r = (r_1, \dots, r_k)$. By definition, $\r$ captures maximal coordinates (with suitable shifts) of shifts in a minimal free resolution of $N$.
By applying a similar line of arguments as that of Theorem~\ref{reg: filter-regular}, we can show that
\begin{align}
 r_l \le \max \{\bfa{l}_N(\x_l), d^\brak{l}(N) \}. \label{app.temp}
\end{align}
Notice that we do not get the equality in (\ref{app.temp}) due to the lack of Nakayama's lemma when $A$ is not necessarily a field.

Now, by virtue of (\ref{app: reg and a-inv}) and (\ref{app.temp}), $\nreg_l(N)$ then gives a bound on the $l$-coordinate (with a suitable shift) in a minimal free resolution of $N$.
\end{remark}

\begin{remark} Let $I$ be a multigraded homogeneous ideal in $T$. Techniques in Section~\ref{s.asympreg} can now be applied, together with (\ref{cohend}), to show that $\nReg(I^nN)$ again is asymptotically a
linear function.
\end{remark}


\end{document}